\documentclass[aps]{revtex4-2}
\usepackage{graphicx}  
\usepackage{dcolumn}   
\usepackage{bm}        
\usepackage{amssymb}   
\usepackage{blindtext}
\usepackage{soul}
\usepackage{float}
\usepackage[textwidth=1.6\marginparwidth,backgroundcolor=yellow,linecolor=orange,textsize=scriptsize]{todonotes}
\usepackage{amsmath}

\usepackage{xspace}
\usepackage{natbib}
\usepackage{subcaption}
\usepackage{ulem}

\newcommand{\mads}{{\sf Mads}\xspace}
\newcommand{\nomad}{{\sf NOMAD}\xspace}
\newcommand{\pynomad}{{\sf PyNOMAD}\xspace}

\def\R{{\mathbb{R}}}

\begin{document}
\boldmath

\title{Blackbox optimization for origami-inspired bistable structures}
\author{Luca Boisneault$^{a,c, \dagger}$, Charles Audet$^{b,d}$, and David Melancon$^{a,c,\dagger}$}
\affiliation{\\a.~Laboratory for Multiscale Mechanics (LM2), Department of Mechanical Engineering, Polytechnique Montréal, 2500 Chemin de Polytechnique Montréal, Québec, H3T 1J4, Canada\\
\\
b.~Department of Mathematical and Industrial Engineering, Polytechnique Montréal, 2500 Chemin de Polytechnique, Montréal, Québec, H3T 1J4, Canada\\
\\
c.~Research Center for High Performance Polymer and Composite Systems (CREPEC), Department of Mechanical Engineering, McGill University, 817 rue Sherbrooke Ouest, Montréal, Québec, H3A 0C3, Canada\\
\\
d.~GERAD, Montréal, Québec, H3T 1J4, Canada\\
\\
$\dagger$ Corresponding authors luca.boisneault@polymtl.ca and  david.melancon@polymtl.ca}

\date{\today}

\begin{abstract}
Bistable mechanical systems exhibit two stable configurations where the elastic energy is locally minimized. 
To realize such systems, origami techniques have been proposed as a versatile platform to design deployable structures with both compact and functional stable states.
Conceptually, a bistable origami motif is composed of two-dimensional surfaces connected by one-dimensional fold lines. 
This leads to stable configurations exhibiting zero-energy local minima. 
Physically, origami-inspired structures are three-dimensional, comprising facets and hinges fabricated in a distinct stable state where residual stresses are minimized. 
This leads to the dominance of one stable state over the other. 
To improve mechanical performance, one can solve the constrained optimization problem of maximizing the bistability of origami structures, defined as the amount of elastic energy required to switch between stable states, while ensuring materials used for the facets and hinges remain within their elastic regime. 
In this study, the Mesh Adaptive Direct Search (\mads) algorithm, a blackbox optimization technique, is used to solve the constrained optimization problem. 
The bistable waterbomb-base origami motif is selected as a case-study to present the methodology. 
The elastic energy of this origami pattern under deployment is calculated via Finite Element simulations which serve as the blackbox in the \mads optimization loop.
To validate the results, optimized waterbomb-base geometries are built via Fused Filament Fabrication and their response under loading is characterized experimentally on a Uniaxial Test Machine. 
Ultimately, our method offers a general framework for optimizing bistability in mechanical systems, presenting opportunities for advancement across various engineering applications.

\end{abstract}

\maketitle

\newpage


\section{Introduction}


Initially an artistic technique of folding paper, origami is now used in engineering to develop deployable systems whose kinematics are embedded directly in the crease pattern \cite{misseroni_origami_2024}.
This has led to the design of sub-millimeter scale mechanical metamaterials \cite{iniguez-rabago_rigid_2022, jamalimehr_rigidly_2022, liu_micrometer-sized_2021, ze_spinning-enabled_2022} and robots \cite{li_fluid-driven_2017, rus_design_2018}  capable of shape reconfiguration as well as meter scale structures \cite{melancon_inflatable_2022, zirbel_hanaflex_2015,zhu_large-scale_2024} that deploy using simple actuation methods.

In engineering, origami research is mainly divided into two categories: rigid foldable \cite{lang_rigidly_2020, feng_designs_2020, mcinerney_discrete_2022} and deformable \cite{martinez_elastomeric_2012, hanna_waterbomb_2014, hanson_controlling_2024}.
Whereas rigid origami can be studied purely from a mechanism point of view, i.e.,~by solving the  equations of motion of rigid bodies \cite{lang_rigidly_2018}, deformable origami requires taking into account the storage of elastic energy to predict deployment.
Mathematically, the folding of deformable origami structures can either be modeled using simple, discrete elements such as bars along a fold line and torsional springs across it \cite{zhu_sequentially_2021}, or using more accurate, finite elements such as thin shells \cite{zhu_review_2022}. While using the Finite Element Method (FEM) to simulate folding provides a rich description of the stored elastic energy inside the origami structure, it comes with an increase in computational cost.

During deployment, deformable origami structures store elastic energy mostly through folding of the hinges and bending of the facets. While hinging energy is typically monotonic, bending energy can be non-monotonic in some origami patterns, leading to  multistability \cite{brunck_elastic_2016}.
This property is defined as the coexistence of two or more equilibrium states where the elastic energy is locally minimized. Recent works have shown that these stable configurations can be accessed via an imposed displacement \cite{hanna_waterbomb_2014, dalaq_experimentally-validated_2022}, magnetic field \cite{novelino_untethered_2020, fang_magneto-origami_2019}, internal pressure \cite{melancon_inflatable_2022, melancon_multistable_2021, zhang_tunable_2023}, or through stimuli-responsive materials \cite{fang_origami-based_2017, zhou_sma-origami_2024}.
Because multistable structures embed self-locking, they offer an advantage over other deployable systems relying on external mechanisms, such as contact \cite{brown_dual-purpose_2022, jamalimehr_rigidly_2022} and spring-loaded devices \cite{holland_development_2016}.
Most of the current multistable origami literature focuses on characterizing the influence of the pattern geometry on multistability \cite{zhang_origami_2018, gillman_design_2019}. 
Some works have put forward optimization as a way to increase bistability, but they are limited to simple beam-based structures \cite{liu_machine_2020} or only maximize geometrical incompatibility \cite{lee_designing_2023, melancon_multistable_2021}. 
However, to transition toward engineering applications, manufacturing parameters, such as panel thickness and hinge type, become important.

This work puts forward a general framework to optimize and take into account multistability when designing origami-inspired structures.
In Section~\ref{sec-Methodology}, the bistable waterbomb pattern \cite{hanna_waterbomb_2014} is chosen as a case study and a modeling representation based on compliant crease is presented.
Its deployment and the associated bistability performance are computed via FEM and validated on 3D printed samples.
The geometry is then parameterized, and the selected design variables are shown to have an impact on the bistable behavior of the origami structure.
The question of finding the best possible geometry is posed as a mathematical optimization problem, in which the objective function consists in maximizing the energy required to switch back from the second to the first state, and is constrained both by the fabrication limitation and the mechanical stress experienced during the deployment phase. 
Evaluating the objective and constraint functions requires launching a time-consuming FEM simulation, which often fails to compute due to instabilities and nonlinearities in the mathematical formulation.  
The resulting optimization problem is solved by the {\sf Nomad} \cite{nomad4paper} implementation of the Mesh Adaptive Direct Search ({\sf Mads}) \cite{AuDe06b} derivative-free constrained blackbox optimization algorithm.
A coupled blackbox-FEM framework is developed to optimize the parameterized model, while taking into account the multiple failed evaluations.
Finally, in Section~\ref{sec-Results}, the optimization process is applied with and without considering manufacturing limits and the resulting geometries are presented. 

\section{Methodology}
\label{sec-Methodology}

\subsection{The origami waterbomb base pattern: a simple bistable structure}

\begin{figure}[h]
    \centering
    \includegraphics[width=\textwidth]{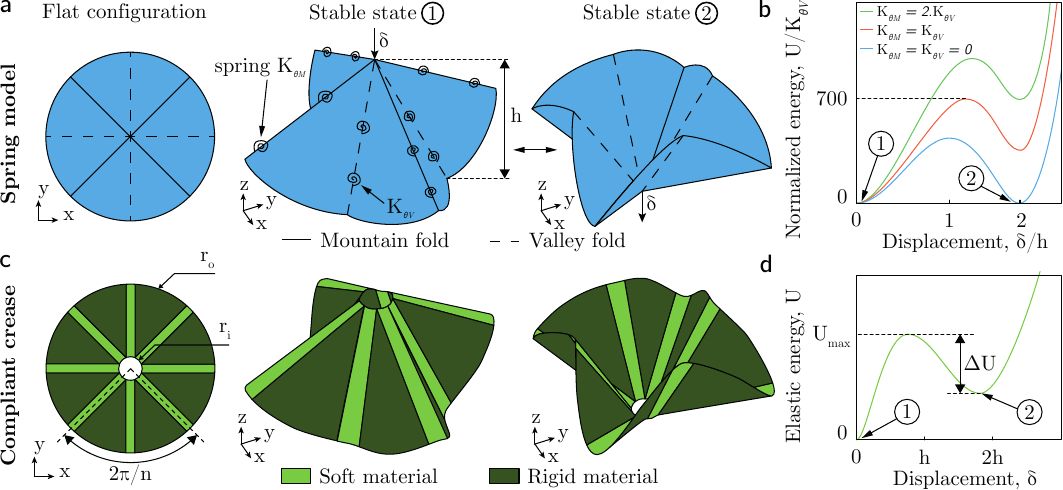}
    \caption{\textbf{Modeling the four-fold waterbomb origami pattern with different levels of complexity.} \textbf{a.}~Simplified representation of the waterbomb with fold lines modeled as torsion springs shown in the flat configuration as well as in both stable states. \textbf{b.}~Energy-displacement curve of the spring model and the effect of increasing the hinging energy by adding springs of stiffness $K_{\theta M}$ and $K_{\theta V}$ at the mountain and valley folds, respectively. \textbf{c.}~Compliant crease representation of the waterbomb with the fold lines modeled as regions of soft material shown in the flat configuration as well as in both stable states. \textbf{d.}~Energy-displacement curve of the compliant crease model.}
    \label{Fig1}
\end{figure}

The waterbomb base pattern is selected as a case study to describe the methodology developed herein to improve the mechanical performance of bistable origami structures. 
This choice is motivated by its simple geometry, ease of fabrication, and extensive researches on its kinematics \cite{imada_geometry_2022}, bistability \cite{hanna_waterbomb_2014}, and potential applications to tune acoustic waves \cite{benouhiba_origami-based_2021}, create logic gates \cite{treml_origami_2018}, build mechanical metamaterials \cite{bai_three-dimensional_2023}, and develop innovative origami-based robots \cite{fang_origami-based_2017}. 
The geometrical description of the waterbomb base is presented in Fig.~\ref{Fig1}. The classical approach to obtain the waterbomb fold is by dividing evenly a flat, circular surface with mountain and valley folds around its geometrical center. 
This axisymmetric pattern usually involves $n=4$ repetitions so that opposite folds are of the same type, mountain or valley, making the structure easier to fold in a cone-like shape of height $h$ (see Fig.~\ref{Fig1}a). When the central node, i.e., the tip of the cone, is pulled down by a distance $\delta$, the structure starts deforming elastically through the bending of the triangular faces and stretching of the fold lines. 
When the $\delta=h$, this stored elastic energy $U$ reaches a maximum.
Passed this point, i.e., for $\delta>h$, the energy $U$ decreases towards a second local minimum.
If the hinges connecting the panels are ideal pivot connections, this second stable state is $z$-symmetrical to the initial configuration so that there is no bending of the triangular faces. 
The bistable behavior of the waterbomb pattern can be characterized mechanically by plotting $U$ as a function of $\delta$, as shown in Fig.~\ref{Fig1}b.
For the case where there is no energy cost associated with rotating the faces along a fold line, face bending prevails and the energy curve shows two stable states with zero energy.
Instead, if torsion springs of stiffness $K_\theta$ are added to model the hinging energy, the second equilibrium state has residual stresses, resulting in a nonzero energy local minimum.
Hanna \textit{et al.} \cite{hanna_waterbomb_2014} have shown that a different stiffness in the mountain and the valley hinges, i.e., $K_{\theta M}$ and $K_{\theta V}$, respectively, will affect the amount of energy required to switch back from the second state to the first one. 
In Fig.~\ref{Fig1}b, their results are reproduced numerically for three different scenarios, i.e., $K_{\theta M}=K_{\theta V}=0$ (blue curve), $K_{\theta M}=K_{\theta V}$, with $K_{\theta V}\neq0$ (red curve), and $K_{\theta M}=2K_{\theta V}$, with $K_{\theta V}\neq0$ (green curve).

In the most simplified representation of the waterbomb pattern, mountain and valley folds are modeled as spring-loaded hinges that connect flat panels. 
A more continuous way of modeling this origami pattern is to represent the folding geometry using compliant creases \cite{zhu_review_2022}.
In this technique, fold lines are substituted with wider and softer regions to allow rotation and the central node is replaced by a hole (see Fig.~\ref{Fig1}c where the dark and light shades correspond to faces and compliant creases, respectively).
This allows to take the hinge width into account, and a higher-order of geometric continuity is implemented, i.e., smooth folds \cite{peraza_hernandez_kinematics_2016} between the faces. 
Applying the compliant crease origami modeling produces the same downside effect as adding torsion springs on the simplified model : the second stable state has nonzero elastic energy. This type of bistable energy curve has two characteristic features: a local maximum of elastic energy between the two stable states, $U_{\max}$, and an energy well depth of the second stable state, $\Delta U$ (see Fig.~\ref{Fig1}d). Here, their ratio is used to quantify the bistability, $\phi$, of the structure:
\begin{equation}
    \phi=\frac{\Delta U}{U_{\max}}.
\end{equation}
When $\phi \rightarrow 0$ the structure becomes marginally bistable. Instead, when $\phi= 1$, the two stable states have the same amount of stored elastic energy.
\unboldmath
\boldmath
\subsection{Simulating bistable origami via the Finite Element Method}  

\begin{figure}[h]
    \centering
    \includegraphics[width=\textwidth]{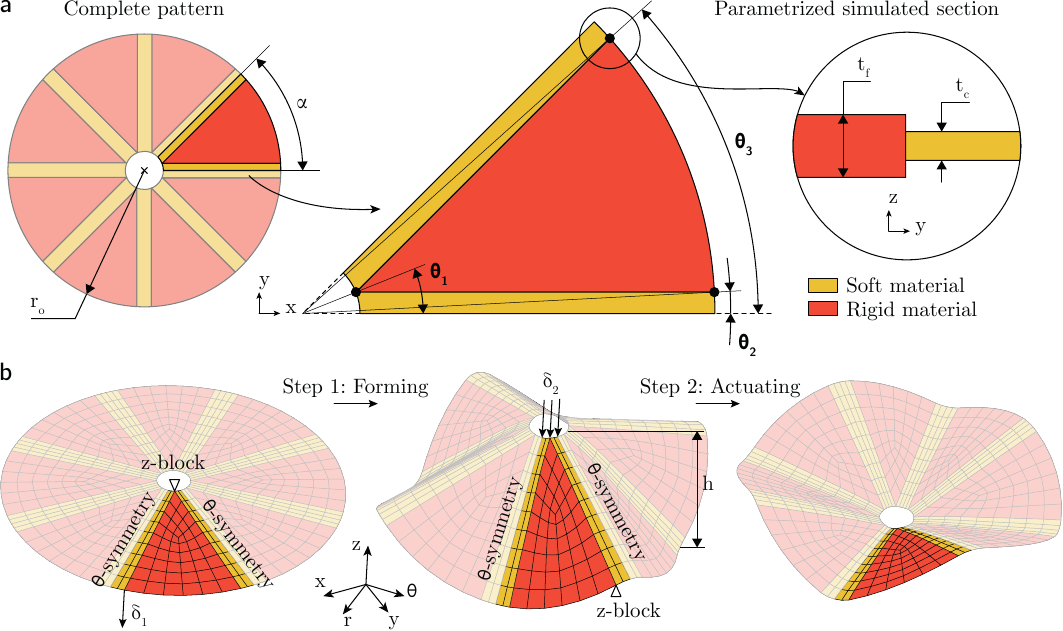}
    \caption{\textbf{Parametrized FEM of the waterbomb origami pattern.} \textbf{a.}~Representation of the geometrical variables $\theta_1$, $\theta_2$, $\theta_3$, $\omega=t_c/t_f$, and $h/r_o$ on the simulated part of the structure. \textbf{b.}~Meshed model of the waterbomb at the beginning and the end of the \textbf{Forming} and the \textbf{Actuation} steps as well as the associated boundary conditions imposed.}
    \label{Fig2}
\end{figure}

In this work, FEM is used to compute the bistability, $\phi$, of the compliant crease waterbomb. The origami pattern is discretized with four-node, linear shell elements (element code S4 in Abaqus Standard 2022) with linear elastic material model with elastic moduli $E_f$, $E_c$, Poisson's ratios, $\nu_f$, $\nu_{c}$, densities, $\rho_f$, $\rho_c$, and elastic limit, $S_Y^f$, $S_Y^c$ for the faces and compliant creases, respectively.
As shown in Fig.~\ref{Fig2}a, for given face and crease materials, five design variables are selected to generate a wide range of geometry for the waterbomb model: three angles, \boldmath$\theta_i$, with $i \in \{1,2,3\}$, shaping the compliant crease, $\omega = t_c/t_f$, the ratio of out-of-plane crease thickness over face thickness, and $h/r_o$, the height of the waterbomb in its first stable state normalized by the outer radius. 
Together, these values constitute the input vector $x=\left(\theta_1/\alpha, \theta_2/\alpha, \theta_3/\alpha, \omega, h/r_o\right) \in \R^5$, with $\alpha = \pi/n$.
The inner radius, $r_i/r_o = 1/6$, and the number of cyclic symmetry, $n = 4$, are fixed to reduce the dimensionality of the design space.
To speed up the computation, only $1/2n$ of the complete pattern is modeled and cyclic boundary conditions are applied on the outer edges. In a cylindrical framework, this means to the two lateral edges cannot move along the $\theta$-axis, as well as rotate around the $r$-axis and the $z$-axis. The FEM simulation is divided into two steps (see Fig.~\ref{Fig2}b):
\begin{itemize}
    \item \textbf{Step-1: Forming.} The waterbomb pattern is deformed from the flat configuration to the deployed state defined by $h/r_o$. 
    To do so, the nodes located on the inner hole are pulled up by a distance $\delta_1 = h$, and the node defined by $\theta_2$, i.e., the node located at the frontier between the face and the crease, is fixed with respect to the $z$-axis.
    The geometry obtained at the end of this step is retrieved and taken without any mechanical stress as the base geometry for the second step.
    \item \textbf{Step-2: Actuation.} The waterbomb pattern is actuated from the first to the second stable state. To do so, the nodes on the inner hole are pulled down by a distance $\delta_2 = 2h$, while the node defined by $\theta_3$, is locked relative to the $z$ translation,
\end{itemize}

During the FEM simulation, the stored elastic energy is obtained by integrating the reaction force with respect to the applied displacement on the nodes of the inner holes during the \textbf{Actuation} step.
In addition, the maximum von Mises stress developed in the structure, $\sigma_{\max}$, as well as its location and associated displacement $\delta$ are extracted from the FEM simulation to ensure the materials remain in their elastic regime, i.e., $\sigma_{\max}<S_Y^i$, with $i\in\{f,c\}$ if the $\sigma_{\max}$ is developed in the faces or crease, respectively.

\begin{figure}[h]
    \centering
    \includegraphics[width=\textwidth]{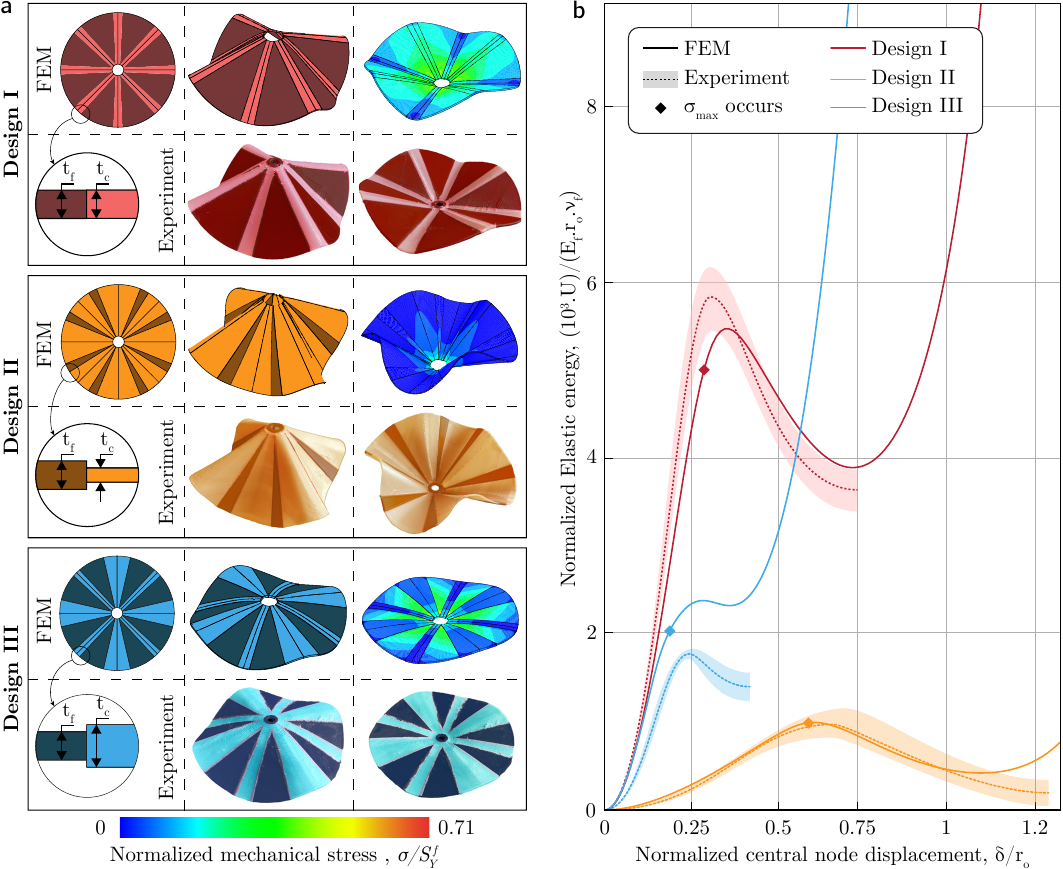}
    \caption{\textbf{Impact of the design variables on the bistability performance of the waterbomb pattern.} \textbf{a.}~The three designs obtained through the parametrization and their geometrical representation in both equilibrium states for the FEM and the 3D printed models. \textbf{b.}~Evolution of the elastic energy during deployment for the three designs showing the comparison between FEM simulations and experimental tests.}
    \label{Fig3}
\end{figure}

To highlight the effect of the geometrical parameters on the bistability of the waterbomb motif, the FEM simulation is conducted on three different patterns: 
\begin{itemize}
    \item \textbf{Design I} with $x_I=\left(\theta_1^I/\alpha, \theta_2^I/\alpha, \theta_3^I/\alpha, \omega^I, h^I/r_o\right)=\left(0.1, 0.5, 0.9, 1.0, 0.6\right)$,
    \item \textbf{Design II} with $x_{II}=\left(\theta_1^{II}/\alpha, \theta_2^{II}/\alpha, \theta_3^{II}/\alpha, \omega^{II}, h^{II}/r_o\right)=\left(0.5, 0.6, 0.7, 0.5, 0.704\right)$,
    \item \textbf{Design III} with $x_{III}=\left(\theta_1^{III}/\alpha, \theta_2^{III}/\alpha, \theta_3^{III}/\alpha, \omega^{III}, h^{III}/r_o\right)=\left(0.31, 0.46, 0.9, 1.5, 0.374\right)$.
    
\end{itemize}
These input vectors are chosen to represent the wide range of feasible geometries. \textbf{Design I} is the standard waterbomb with parallel creases and equal thickness between facets and creases. Differently, \textbf{Design II} includes creases wider, but thinner than facets. Finally, \textbf{Design III} alternates narrow and wide creases which are thicker than the facets.
For each design, the top and front views as well as the two stable configurations and the von Mises stress field in the second stable state are shown in Fig.~\ref{Fig3}a.

Here, the ratios $E_f/E_c = 21.67$ and $\nu_f/\nu_c = 0.78$ are considered in the numerical simulations.
The evolution of the elastic energy during the deployment of each design is presented in Fig.~\ref{Fig3}b and reveals that changes in the size and shape of the compliant creases can affect drastically the bistable performance of the waterbomb pattern.
The elastic energy is normalized with respect with the outer radius $r_o$ of the pattern, as well as the facet's material properties $E_f$ and $\nu_f$.
For \textbf{Design I}, the bistable performance is characterized by $\phi^I = 28.85\%$, with $\Delta U^{I}/(E_fr_o\nu_f)=1.56\times 10^{-3}$ and $U^{I}_{\max}/(E_fr_o\nu_f)=5.41\times 10^{-3}$. 
\textbf{Design II} exhibits an increase of bistable performance with $\phi^{II} = 57.74\%$, $\Delta U^{II}/(E_fr_o\nu_f)=0.57\times 10^{-3}$, and $U^{II}_{\max}/(E_fr_o\nu_f)=0.98\times 10^{-3}$. Finally, \textbf{Design III} displays marginal bistability with $\phi^{III} = 2.4\%$, $\Delta U^{III}/(E_fr_o\nu_f)=0.05\times 10^{-3}$, and $U^{III}_{\max}/(E_fr_o\nu_f)=2.36\times 10^{-3}$.
Note that the geometry of the creases affect not only the multistability ratio $\phi$, but also the maximum elastic energy $U_{max}$, the barrier of energy in the second stable state $\Delta U$, and the displacement $\delta$ required to switch to the second stable state. 

Additionally, for the three design, the maximum mechanical stress is developed right before the local maximum of energy (see diamond markers in Fig.~\ref{Fig3}b) and is located near the hole and on the stiff faces. This maximal value, normalized by the elastic limit of the facets material, is $\sigma^{I}_{\max}/S_Y^f=0.708$ MPa, $\sigma^{II}_{\max}/S_Y^f=0.685$ MPa and $\sigma^{III}_{\max}/S_Y^f=0.473$ MPa for the three designs.
The location of the maximum stress for the second stable state stays the same (see the contour maps in Fig.\ref{Fig3}a), but one notes that a higher mechanical stress is associated with a lower bistability performance : $\sigma^{I}_{state2}/S_Y^f=0.304$ MPa, $\sigma^{II}_{state2}/S_Y^f=0.266$ MPa and $\sigma^{III}_{state2}/S_Y^f=0.432$ MPa.

\subsection{Experimental validation}

To validate the FEM simulations, physical prototypes of \textbf{Designs I-III} are fabricated using the Fused Filament Fabrication (FFF) method.
FFF enables fast prototyping of multi-material origami patterns \cite{ye_multimaterial_2023}.
Here, the crease regions are printed with thermoplastic polyurethane (TPU from Eryone with $E_c=120$ MPa, $\nu_c=0.45$, $\rho_c=1200$ kg/m\textsuperscript{3} and $S_Y^c = 50$ MPa) and the faces with polylactic acid (PLA from Raise3D with $E_f=2600$, $\nu_f=0.35$, $\rho_f=1040$ kg/m\textsuperscript{3} and $S_Y^f = 50$ MPa).
Using the initial deformed shape obtained from the FEM simulation of the waterbomb (\textbf{Step-1: Forming}), a CAD model is generated. 
Interlocking geometry are added to improve bonding between the soft and the rigid regions \cite{kuipers_itil_2022}.
This way, when the model is sliced, the interface between faces and creases will be composed of layers of an alternating set of PLA and TPU layers (see the Supplementary Materials section S1 for additional details on the fabrication of the physical prototypes).

To experimentally measure the stored elastic energy during folding, the physical prototypes are tested on a uniaxial test machine (MTS Insight Electromechanical 50), fitted with a $100$ N load cell.
Inspired by previous experimental work on the waterbomb pattern \cite{hanna_waterbomb_2014}, each specimen is placed on custom-build supports that emulate the same conditions as the numerical simulation, i.e., it is installed on triangular guiding rails preventing it from rotating around the $z$-axis, and the central hole is attached to the load cell through a fixed bolt to measure the reaction force during the deployment (see the Supplementary Materials S2 for additional details on the experimental testing).
The crosshead imposes a vertical displacement of $\delta_{exp}$ at a rate of $0.05$ mm/s, to create quasi-static conditions, until the second stable state is reached.
After the test, the reaction force $F_{\exp}$ is integrated along the displacement to obtain the experimental elastic energy $U_{\exp}$ of the waterbomb prototype:
\begin{equation}\label{eq:continuous}
    U_{\exp}(\delta_{exp})=\int_{0}^{\delta_{exp}} F_{\exp}(\delta) \, \mathrm{d}\delta.
\end{equation}
In Fig.~\ref{Fig3}a, representative specimens of each prototype are shown in both stable states, and their measured energy-displacement curves are plotted as doted line in Fig.~\ref{Fig3}b with the standard deviation from five tests specimens shown as shaded areas.
From the comparison with the simulated model, one can assess that the deformed shapes obtained with the 3D printed sample match qualitatively the computation both in the first and the second stable states (see Fig.~\ref{Fig3}a). Quantitatively, \textbf{Design I} shows the closest match between predicted and measured energy landscapes with a relative error of $\epsilon_{U_{\max}}^I=5.5\%$ on the maximal amount of energy $U_{\max}$ and an error on the displacement of the central node in the second stable state $\epsilon_{\delta/r_o}^I=4.1\%$.
\textbf{Design II} also shows good agreement between simulations and experiments with $\epsilon_{U_{\max}}^{II}=2.9\%$ and $\epsilon_{\delta/r_o}^{II}=18.4\%$.
However, there are discrepancies between the predicted and measured energy landscapes for \textbf{Design III} with $\epsilon_{\delta/r_o}^{III}=15.9\%$ and $\epsilon_{U_{\max}}^{III}=25.7\%$. This deviation could be attributed to the boundary condition imposed during the testing phase which can slightly differ from the one set numerically.
In fact, the larger soft creases of \textbf{Design III} close to the hole where the screw is fixed for mechanical testing could add additional compliance that is not modeled in the FEM simulations. 



\subsection{Optimizing the mechanical performance of bistable origami}

Optimization algorithms can be used to tune the geometrical parameters of the waterbomb origami with compliant creases in order to address the loss of bistable performance seen in Fig.~\ref{Fig3}b.
Consider a simulation that takes an input vector $x$ containing the five design variables, and outputs the bistable performance of the associated structure, $\phi$, as well as the maximum von Mises stress experienced by the structure, $\sigma_{\max}$.
The input vector $x \in \mathbb{R}^5$ is bounded by the two vectors $l_b$ and $u_b$, respectively the lower and the upper boundaries.
This ensures that the optimization will not diverge and deliver unrealistic results.
To avoid mechanical failure, the maximum stress experienced by the structure $\sigma_{\max}$ must not exceed the yield strength, $S_Y$, of the material.
The corresponding constrained optimization problem is then formulated as :
\begin{equation}
\begin{aligned}
\max_{x \in \mathbb{R}^5} \quad & \phi=\frac{\Delta U}{U_{\max}}\\
\textrm{s.t.} \quad & \sigma_{\max}\leq\ S_Y\\
  &  l_b \leq x \leq u_b.  \\
\end{aligned}
\end{equation}

This optimization problem may be regarded as a blackbox : {at each iteration $k$}
    only the input $x_k$ and output $\phi^k$ and $\sigma_{\max}^k$ data are known, and the time-consuming FEM simulations are considered hidden.
In the present case, derivatives are difficult to obtain due to the numerous fails in the computation, therefore derivative-free optimization techniques~\cite{AuHa2017} are required.

The \nomad blackbox optimization software, which implements the \mads algorithm~\cite{AuDe2006}, is chosen to maximize the objective function while taking into account the constraints of the problem.
Thanks to the dynamic adaption of the size of the searching space between each iteration, \mads allows an efficient exploration of the design space.
Additionally, \nomad has proven successful in the case of blackbox with a long computational time \cite{AuOr06a} and with large part of the design space covered by hidden constraints, i.e.,~sets of points that cannot be computed or that do not output numerical values when fed to the blackbox \cite{ChKe2016}.
Numerous studies in fields such as biomedical, aerospace or electrical engineering, have successfully used \nomad to solve optimization problems \cite{AlAuGhKoLed2020}.
This software is also provided with a python package, \pynomad, allowing easy communication with the FE software.

\begin{figure}[h!]
    \centering
    \includegraphics[width=\textwidth]{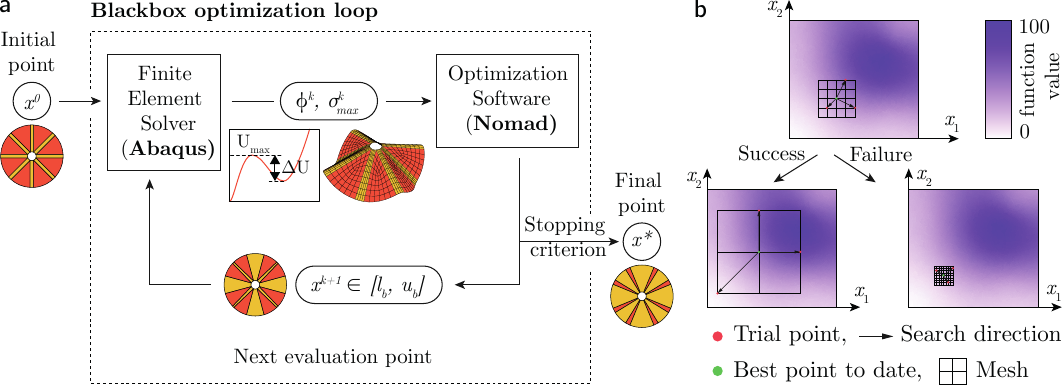}
    \caption{\textbf{Strategy to optimize the bistability of origami-inspired structures.} \textbf{a.}~Coupling of the FEM blackbox with the optimization algorithm \nomad. \textbf{b.}~Successive polling steps on an arbitrary 2D function with \mads depending if the previous iteration was a success or a failure.}
    \label{Fig4}
\end{figure}

In the present problem, each call to the blackbox requires an average of $110$ seconds to compute when successful, and computation fails on $30$\% of the calls, which makes \nomad a suitable solution to solve the problem. 
Computations are made on an Intel Core i9-9900K processor.
Here, every evaluation requires a different computation time, depending on how well the FE solver performs on the model defined by a given input vector.
The optimization process is initiated with \textbf{Design I}, a geometry with parallel creases and uniform thickness across the structure.
This geometry is set in \nomad as the initial point $x^0 = \left(\theta_1^0/\alpha, \theta_2^0/\alpha, \theta_3^0/\alpha, \mu^0, h^0/r_o\right) = \left(0.1, 0.5, 0.9, 1.0, 0.6\right)$, and is known to provide a bistable performance of $\phi^I=28.85\%$.
\nomad sends these parameters to the blackbox for the first evaluation and the FE software, which computes the stored elastic energy and von Mises stress during folding, outputs back to \nomad both the bistable performance $\phi$ as well as the maximum value of stress $\sigma_{\max}$ (Fig.~\ref{Fig4}a).
The next evaluation points are determined by selecting  $N+1$ random points (with $N$ is the dimension of the design space) on a grid centered on the best evaluation yet.
A complete iteration of the optimization algorithm consists of the evaluation of these $N+1$ points.
If the multistable performance $\phi$ of the evaluation $k$ is better than any of the previous best evaluations, $x^k$ and the corresponding $\phi^k$ become the new champion and the size of the grid that determines the next iteration points is expanded.
However, if the computation does not result in an improvement of the objective function, the size of the grid is reduced for the next iteration. 
If the mechanical stress exceeds the yield stress limit, the point is discarded and cannot be the final output value of the optimization, i.e., it is considered a failed evaluation.
The algorithm continues until it reaches a maximum of $1000$ evaluations.
The optimization process is schematized in Fig.~\ref{Fig4}b.

\newpage
\section{Optimization results}
\label{sec-Results}

\begin{figure}[h]
    \centering
    \includegraphics[width=\textwidth]{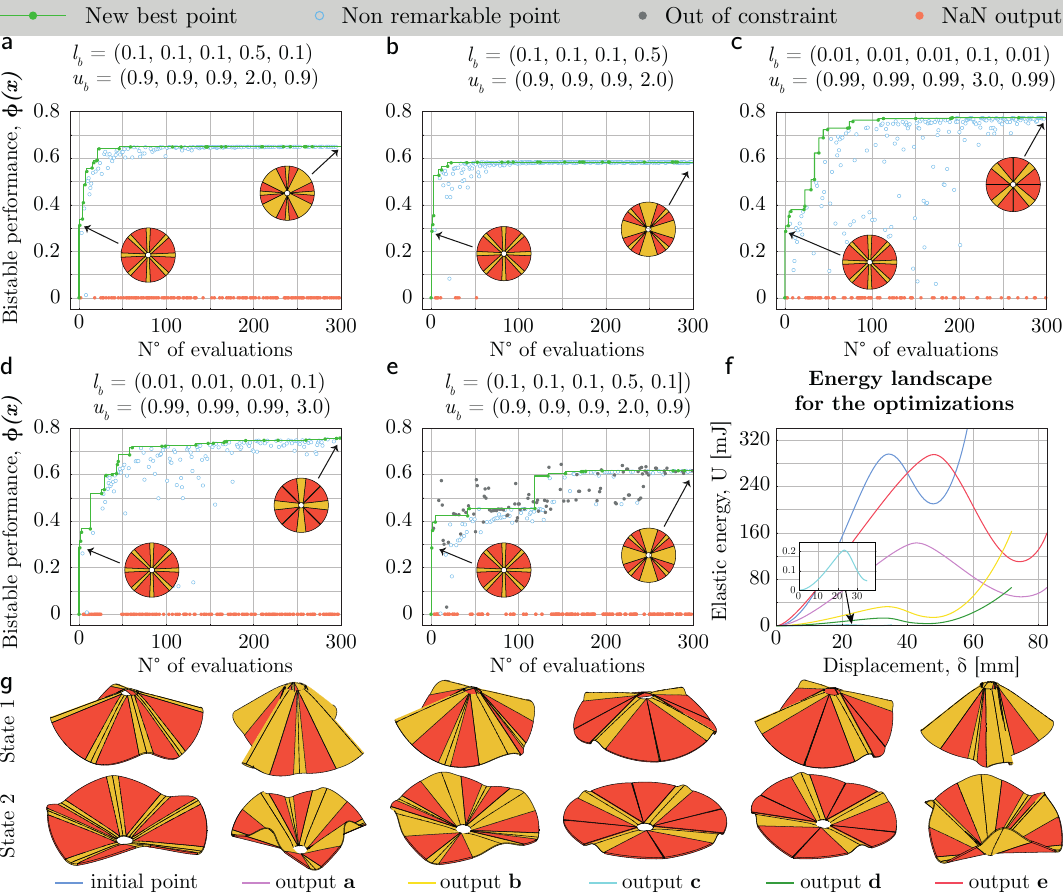}
    \caption{\textbf{Results of the blackbox optimization.} Bistable performance $\phi(x)$ as a function of the number of evaluations when taking into account manufacturing limits while leaving free \textbf{(a)} and fixing \textbf{(b)} the initial height $h/r_o$. Effect of increasing the range of the design variables on $\phi(x)$ while leaving free \textbf{(c)} and fixing \textbf{(d)} the initial height $h/r_o$. \textbf{e.}~Effect of adding an additional constraint on the maximum energy developed during deployment $U_{\max}$. \textbf{f.}~Evolution of the elastic energy during deployment for the initial geometry and the five optimized geometries \textbf{(a)}-\textbf{(e)}. \textbf{g.}~The two stable states of the waterbomb for the initial geometry and the five optimized geometries.}
    \label{Fig5}
\end{figure}

The optimization is launched considering manufacturing limits, i.e., ensuring a minimum hinge width of $0.4$ mm, the diameter of the nozzle used for 3D printing the samples. 
The associated mathematical lower and upper boundaries $l_b$ and $u_b$ as well as the evolution of the bistable performance $\phi(x)$ with the number of evaluations are displayed in Fig.~\ref{Fig5}a.
The design variables of the optimize geometry are $x^a = \left(\theta_1^a/\alpha, \theta_2^a/\alpha, \theta_3^a/\alpha, \omega^a, h^a/r_o\right) = \left(0.1000, 0.1000, 0.5421, 0.5296, 0.7359\right)$ with an associated bistable performance of $\phi_a = 64.8\%$.
The performance here is more than doubled if compared from the initial design (see the energy landscape and stable configurations of the initial and optimized geometries in Figs.~\ref{Fig5}f-g). 
Around $30\%$ of the evaluations performed ended up failing.
This result can be linked to the increase of the last variable, $h/r_o$. Tall and narrow waterbomb folds are associated with high geometrical frustration during reconfiguration and this can introduce high nonlinearities in the numerical simulations.
To reduce the number of failed evaluations, a second optimization is launched with the initial height fixed to $h/r_o=0.6$ along the process.
This optimization produces the final vector $x^b = \left(0.4366, 0.9000, 0.9000, 0.5000, 0.6000\right)$ and the bistable performance $\phi_b=58.4\%$ as shown in the convergence plot of Fig.~\ref{Fig5}b.
While this represents a loss of $6.4\%$ compared to the results in Fig.~\ref{Fig5}a, it still shows a two-fold increase with respect to the initial design.
In addition, the number of failed computations goes from $30\%$ to only $1\%$ over $1000$ evaluations. For the two different optimizations, the convergence is fast with $99\%$ of the final performance already reached after only $25$ evaluations.

From the insets showing the final geometries in Figs.~\ref{Fig5}a-b and the stable states in Fig.~\ref{Fig5}g, one notes that the optimization leads to configurations which reach the mathematical constraints for certain angles $\theta_i/\alpha$, i.e., narrow compliant creases. To investigate the potential gain associated with a manufacturing technique with higher resolution, the lower/upper boundaries on $\theta_i/\alpha$ are decreased/increased and two additional optimizations are launched: one with all five design variables (Fig.~\ref{Fig5}c) and one with the initial height fixed to $0.6$ (Fig.~\ref{Fig5}d). For both cases, the valley folds of the final geometry become even narrower to increase the bistable performance to $\phi_c = 77.7\%$ and $\phi_d = 76.8\%$ (see the corresponding convergence graphs in Figs.~\ref{Fig5}c-d, energy curves in Fig.~\ref{Fig5}f, and stable configurations in Fig.~\ref{Fig5}g). Importantly, for these two cases, the increase in bistability performance is linked to steep lowering in the elastic energy. While for the initial design $U_{\max}^0=36.93$ \textbf{mJ}, the two optimal geometries shown in Figs.~\ref{Fig5}c-d display $U_{\max}^{c}=0.21$ \textbf{mJ} and $U_{\max}^d=1.57$ \textbf{mJ}, respectively. For load-bearing applications, one may want to design bistable origami structures which maximize bistability while being able to develop a high amount of elastic energy during deployment. 
To do so, two methods can be applied : using stiffer materials to manufacture the structure or adding $U_{\max}$ a new mathematical constraint in the optimization.
Here, the second approach is implemented with the added constraint forcing the optimization to seek for designs with at least the same $U_{\max}$ as the initial design :
\begin{equation}
    U_{\max} \geq U_{\max}^0
\end{equation}
In Fig.~\ref{Fig5}e, the last case is presented leading to $x^e  = \left(0.4425, 0.8921, 0.8999, 0.6976, 0.8690\right)$ with a bistable performance of $\phi_e = 62.3\%$. For this last scenario, the convergence of the optimization displays two successive plateaus, caused by the additional difficulty for the algorithm to find geometries that sustain an acceptable level of maximum elastic energy.
The output geometry is almost identical to the one obtained in Figs~\ref{Fig5}a-b, but with a higher initial height and thicker soft regions, characteristics that affect the order of magnitude of the elastic energy response.

\newpage
\section{Discussion}

In this work, an optimization framework is developed to improve the bistability performance of origami-inspired structures and applied to the waterbomb base pattern. The optimization results highlight a two-fold increase in bistability performance from the classic straight crease waterbomb pattern to a more complex geometry with uneven creases. The methodology developed here is general and can be applied to other bistable origami-inspired structures (see the Supplementary Materials section S3 for more details).

The presented framework is adaptable and could be further improved. First, implementing a bar-and-hinge model \cite{zhu_bar_2020} as surrogate computation model~\cite{BoDeFrSeToTr99a,AuLedSa21} could speed up the optimization. 
Additional variables could be easily introduced in the algorithm, e.g., a categorical variable~\cite{G-2022-11} that would determine the material used for each region of the origami pattern, or curved creases~\cite{FLORES2022101909} to get more flexibility on the crease shape.
As shown with the optimization results in Fig.~\ref{Fig5}c-d, relaxing the optimization bounds, which could be possible with other, high resolution fabrication techniques such as composite laminate \cite{suzuki_origami-inspired_2020}, could further increase bistability. 
Finally, the optimization strategy could be extended to multistable origami structures, i.e., with more than two stable states.
To do so, one could use multi-objective optimization, but this technique tends to lack efficiency and often designers have to prioritize one objective over the other~\cite{AuSaZg2010a}.
As multistable origami structures are often made of an assembly of building blocks, e.g., kresling arrays \cite{wang_multi-triangles_2023}, the optimization could be carried both locally on individual components and globally to ensure geometric compatibility. 

\section*{Acknowledgments}
\small{
\textbf{Funding:} This work was funded by Audet's NSERC Canada Discovery Grant 2020--04448 and by Melancon's research start-up grant from Polytechnique Montréal. \textbf{Author contributions:} L.B., C.A., and D.M. proposed and developed the research idea. L.B. performed the FEM simulations, fabricated the physical prototypes, and developed and ran the blackbox optimization. L.B., C.A. and D.M. wrote the paper. C.A. and D.M. supervised the research. \textbf{Competing interests:} The authors declare no conflict of interest. \textbf{Data and materials availability:} The data that support the findings of this study are available on request from the corresponding author. The optimization code described schematically in Fig.~\ref{Fig4} is available on Github at \texttt{https://github.com/lm2-poly/OriMads}
}

\unboldmath

\bibliographystyle{apalike}
\bibliography{bibliography}

\end{document}